\documentclass[a4paper,10pt]{article}

\usepackage{amssymb, amscd, verbatim, amsthm}
\usepackage[all]{xy}
\usepackage{hyperref}
\usepackage{amsmath}

\begin{document}

\title{Big monodromy theorem for abelian varieties over finitely
generated fields}

\author{Sara Arias-de-Reyna\\
        Institut f{\" u}r Experimentelle\\
        Mathematik,\\
        45326 Essen, Germany\\
        \texttt{\small sara.arias-de-reyna@uni-due.de}
        \vspace{3mm}
\and
        Wojciech Gajda\\
        Department of Mathematics,\\
        Adam Mickiewicz University,\\
        61614 Pozna\'{n}, Poland\\
        \texttt{\small gajda@amu.edu.pl}\\
        \vspace{3mm}
\and
        Sebastian Petersen\footnote{the corresponding author}\\
        Universit\"at Kassel,\\
        Fachbereich Mathematik,\\
        34132 Kassel,
        Germany\\
        \texttt{\small petersen@mathematik.uni-kassel.de}
}
\parindent0em
\parskip1em

\maketitle

\newtheorem{leer}{}[section]
\newtheorem{thm}[leer]{Theorem}
\newtheorem{conj}[leer]{Conjecture}
\newtheorem{defi}[leer]{Definition}
\newtheorem{rema}[leer]{Remark} 
\newtheorem{prop}[leer]{Proposition}
\newtheorem{lemm}[leer]{Lemma}
\newtheorem{coro}[leer]{Corollary}
\newtheorem{quest}[leer]{Question}
\newtheorem{claim}[leer]{Claim}

\newtheorem*{abcSatz}{Main Theorem}
\newtheorem*{abcFolgerung}{Corollary}

\newcommand{\prm}{\ {\mathrm{prime}}}

\newcommand{\ilim}{\mathop{\varinjlim}\limits}
\newcommand{\plim}{\mathop{\varprojlim}\limits}

\newcommand{\OO}{{\cal O}}
\newcommand{\PP}{{\cal P}}
\renewcommand{\AA}{{\cal A}}
\newcommand{\BB}{{\cal B}}
\newcommand{\CC}{{\cal C}}
\newcommand{\DD}{{\cal D}}
\newcommand{\LL}{{\cal L}}
\newcommand{\UU}{{\cal U}}
\newcommand{\TT}{{\cal T}}
\newcommand{\MM}{{\cal M}}
\newcommand{\FF}{{\cal F}}
\newcommand{\NN}{{\cal N}}
\newcommand{\RR}{{\cal R}}
\newcommand{\KK}{{\cal K}}
\newcommand{\XX}{{\cal X}}
\newcommand{\YY}{{\cal Y}}
\renewcommand{\SS}{{\cal S}}
\newcommand{\Mon}{{\cal M}}

\newcommand{\Occ}{{\mathrm{Occ}}}
\newcommand{\sep}{_{\mathrm{sep}}}
\newcommand{\bad}{_{\mathrm{bad}}}
\newcommand{\tor}{_{\mathrm{tor}}}
\newcommand{\ab}{_{\mathrm{ab}}}
\newcommand{\Spec}{\mathrm{Spec}}
\newcommand{\Eig}{\mathrm{Eig}}
\newcommand{\dro}{\mathrm{drop}}
\newcommand{\Sym}{\mathrm{Sym}}
\newcommand{\Sp}{\mathrm{Sp}}
\newcommand{\GL}{\mathrm{GL}}
\newcommand{\amp}{\mathrm{amp}}
\newcommand{\Emb}{\mathrm{Emb}}
\newcommand{\GSp}{\mathrm{GSp}}
\newcommand{\PSp}{\mathrm{PSp}}
\newcommand{\SL}{\mathrm{SL}}
\newcommand{\Hom}{\mathrm{Hom}}
\newcommand{\Alg}{\mathrm{Alg}}
\newcommand{\ggT}{{\mathrm{ggT}}}
\newcommand{\can}{{\mathrm{can}}}
\newcommand{\trdeg}{\mathrm{trdeg}}
\newcommand{\Mor}{\mathrm{Mor}}
\newcommand{\mf}{\mathfrak}
\newcommand{\ol}[1]{\overline{#1}}
\newcommand{\Var}{\mathrm{Var}}
\newcommand{\AbVar}{\mathrm{AbVar}}
\newcommand{\Res}{\mathrm{Res}}
\newcommand{\Grp}{\mathrm{Grp}}
\newcommand{\Sch}{\mathrm{Sch}}
\newcommand{\Inert}{\mathrm{Inert}}
\newcommand{\Iso}{\mathrm{Iso}}
\newcommand{\Cur}{\mathrm{Cur}}
\newcommand{\Alb}{\mathrm{Alb}}
\newcommand{\Rat}{\mathrm{Rat}}
\newcommand{\Aut}{\mathrm{Aut}}
\newcommand{\End}{\mathrm{End}}
\newcommand{\DivCor}{\mathrm{DivCor}}
\newcommand{\degr}{\mathrm{deg}}
\newcommand{\ran}{\mathrm{rg}}
\newcommand{\rk}{\mathrm{rk}}
\newcommand{\ke}{\mathrm{ker}}
\newcommand{\coke}{\mathrm{coker}}
\newcommand{\im}{\mathrm{im}}
\newcommand{\chara}{\mathrm{char}}
\newcommand{\rot}{\mathrm{rot}}
\newcommand{\ord}{{\mathrm{ord}}}
\newcommand{\Rest}{\ \mathrm{Rest}\ }
\newcommand\F{\mathbb{F}}
\newcommand\Gal{\mathrm{Gal}}
\newcommand\pr{\mathrm{pr}}
\newcommand\Id{\mathrm{Id}}
\newcommand\Stab{\mathrm{Stab}}
\newcommand\calB{\mathcal{B}}

\newcommand{\dR}[4]{#1\div #2=#3\ \Rest #4}

\newcommand{\TS}{TS}

\def\Pp{{\mathbb{P}}}
\def\Ff{{\mathbb{F}}}
\def\Rr{{\mathbb{R}}}
\def\Cc{{\mathbb{C}}}
\def\Qq{{\mathbb{Q}}}
\def\Zz{{\mathbb{Z}}}
\def\Aa{{\mathbb{A}}}
\def\Nn{{\mathbb{N}}}
\def\Gg{{\mathbb{G}}}

\begin{abstract}
An abelian variety over a field $K$ is said to have {\em big monodromy}, if the image of the Galois representation on $\ell$-torsion points, for almost all primes $\ell,$ contains the full symplectic group.
We prove that all abelian varieties over a finitely generated field $K$ with the endomorphism ring $\Zz$
and semistable reduction of toric dimension one at a place of the base field $K$ have big monodromy. We make
no assumption on the transcendence degree or on the characteristic of $K$.
This generalizes a recent result of Chris Hall.  
\end{abstract}


\section*{Introduction}
\addcontentsline{toc}{section}{Introduction}
It has been known in number theory, since times immemorial that Galois representation attached 
to the action of the absolute Galois group on torsion points of an abelian group scheme carries 
a lot of basic arithmetic and geometric information. The first aim which one encounters naturally, 
while studying such representations is to determine their images in terms of linear algebraic groups. 
There exists a vast variety of results in the literature concerning computations of Galois 
representations for abelian varieties defined over number fields and their applications to some classical 
questions such as Hodge, Tate and Mumford-Tate conjectures, see for example \cite{serretovigneras}, \cite{murty}, 
\cite{ribet}, \cite{bgk3} or \cite{vasiu}. In this paper we are interested in computing images of Galois representations attached to abelian varieties defined over finitely generated fields in arbitrary characteristic, i.e., to families of abelian varieties.  

Let $K$ be a field and denote by $G_K$ its absolute Galois group. Let $A/K$  
be an abelian variety and $\ell\neq\chara(K)$ a prime number. We denote by 
$\rho_{A[\ell]}:G_K\longrightarrow \Aut(A[\ell])$ the Galois representation attached
to the action of $G_K$ on the $\ell$-torsion points of $A.$ 
We define $\Mon_K(A[\ell]):=\rho_{A[\ell]}(G_K)$ and call
this group {\em the mod-{$\ell$} monodromy group of $A/K$}. 
We fix a polarization and denote by $e_\ell{:}A[\ell]\times A[\ell]\to \mu_\ell$ the corresponding
Weil pairing. Then $\Mon_{K}(A[\ell])$ is a subgroup of the group of symplectic similitudes $\GSp(A[\ell],e_{\ell})$ of the Weil pairing. We will say that $A/K$ {\em has big monodromy} if there exists a constant $\ell_0$ such that
$\Mon_{K}(A[\ell])$ contains the symplectic group $\Sp(A[\ell], e_\ell),$ for every prime number $\ell\geq \ell_0.$
Note that the property of having big monodromy does not depend on the choice of the polarization.

Certainly, the most prominent result on computing monodromy groups is the 
classical theorem of Serre (cf. \cite{serre1984}, \cite{serre1985}): 
{\em If $A$ is an abelian variety over a finitely generated field $K$ of characteristic zero with $\End(A)=\Zz$ and $\dim(A)=2, 6$ or odd, then $A/K$ has big monodromy.} In this paper we consider monodromies for abelian varieties over finitely generated fields which have been recently investigated by Chris Hall \cite{hall2007}, \cite{hall2008}. To simplify notation, we will say that 
an abelian variety $A$ over a finitely generated field $K$ {\em is of Hall type}, if $\End(A)=\Zz$ and $K$ has a discrete valuation at which $A$ has semistable reduction of toric dimension one.

In the special case, when $K=F(t)$ is a rational function field over another
finitely generated field, it has been shown by Hall that certain hyperelliptic Jacobians
have big monodromy; namely the Jacobians $J_C$ of hyperelliptic curves
$C/K$ with affine equation $C: Y^2=(X-t)f(X)$, where $f\in F[X]$ is a monic squarefree
polynomial of even degree $\ge 4$ (cf. \cite[Theorem 5.1]{hall2007}).
Furthermore, Hall has proved recently \cite{hall2008} the following
theorem which in our notation reads: {\em If $K$ is a {\em global} field, then
every abelian variety $A/K$ of Hall type has big monodromy.}
We strengthen these results in our main theorem as follows.

\begin{abcSatz}\label{ThmB}[cf. Thm. \ref{main1}] If $K$ is a
finitely generated field (of arbitrary characteristic) and $A/K$ is
an abelian variety of Hall type, then $A/K$ has big monodromy.
\end{abcSatz}

Our proof of the main Theorem follows Hall's proof of \cite{hall2008} to some extent, e.g., we have borrowed
a group theory result from \cite{hall2008} (cf. Theorem \ref{stronghall}). 
In addition to that we had to apply a substantial quantity of new methods to achieve the extension 
to all finitely generated fields, such as for instance finite generation properties of fundamental groups 
of schemes and Galois theory of certain division fields of abelian varieties, which are gathered in Section 2 and Section 3 of the paper. Furthermore, at a technical point in the case $\chara(K)=0$, we perform a tricky reduction argument (described in detail in Section 3) at a place of $K$ whose residue field is a number field. The paper carries an appendix with a self contained proof of the group theoretical Theorem \ref{stronghall} 
due to Hall, which can be of independent value for the reader. 

In the positive characteristic case, to the best of our knowledge, Theorem A is the first result which provides 
a full description of the monodromy groups for a class of abelian varieties of arbitrary dimension, defined over a finitely generated field of transcedence degree bigger than one. 

Theorem A plays an important role in our paper \cite{agppreprint0}, where we make progress 
on the conjecture of Geyer and Jarden (cf. \cite{geyerjarden1978}) on torsion of abelian varieties over 
large algebraic extensions of finitely generated fields. 

As a further application, we combine our monodromy computation with recent results of 
Ellenberg, Hall and Kowalski in order to obtain the following
result on endomorphism rings and simplicity of fibres in certain families of abelian varieties.
If $K$ is a finitely generated transcendental extension of another field $F$ and $A/K$ is an abelian variety, then we
call $A$ {\em weakly isotrivial with respect to $F$}, if there is
an abelian variety $B/\widetilde{F}$ and an $\widetilde{K}$-isogeny $B_{\widetilde{K}}
\to A_{\widetilde{K}}$.

\begin{abcFolgerung}\label{ThmE} [cf. Cor. \ref{ehk}]
Let $F$ be a finitely generated field and
$K=F(t)$ the function field of $\Pp^1/F$. Let $A/K$ be an abelian
variety. Let $U\subset \Pp^1$ be an open subscheme such that
$A$ extends to an abelian scheme $\AA/U$. For $u\in U(F)$ denote
by $A_u/F$ the corresponding special fiber of $\AA$.
Assume that $A$ is not weakly isotrivial with respect to $F$ and that either of the conditions i) or ii) listed below is satisfied.
\begin{enumerate}
\item[i)] $A$ is of Hall type.
\item[ii)] $\chara(K)=0$, $\End(A)=\Zz$ and $\dim(A)=2, 6$ or odd.
\end{enumerate}
Then the sets:
$$X_1:=\{u\in U(F)\,|\quad\End(A_u)\neq \Zz\}$$
and
$$X_2:= \{u\in U(F)\,|\quad\mbox{$A_u/F$ is not geometrically simple}\}$$
are finite.
\end{abcFolgerung}

Note that Ellenberg, Elsholtz, Hall and Kowalski proved the statement of the Corollary in the
special case when $A$ is the Jacobian variety of the hyperelliptic curve given by the affine equation
$Y^2=(X-t)f(X)$, with $f\in F[X]$ squarefree and monic of even degree
$\ge 4$ (cf. \cite[Theorem 8]{ehk2009}). It is the case, where the monodromy of $A$ is known by
\cite[Theorem 5.1]{hall2007}. We obtain part (i) of the Corollary as a consequence of the main
theorem, our Proposition \ref{geo} below and also Propositions 4 and 7 of \cite{ehk2009}.
In order to prove (ii) we use Serre's Theorem \cite{serre1984}, \cite{serre1985} instead of the main Theorem.

We warmly thank Gerhard Frey, Dieter Geyer, Cornelius Greither and Moshe Jarden for conversations and useful comments on the topic of this paper. The mathematical content of the present work has been much influenced by seminal
results of J.-P. Serre contained in \cite{serre1984}, \cite{serre1985}, \cite{serretoribet}, \cite{serretovigneras} and by the inspiring paper \cite{hall2008} of C.Hall. We acknowledge this with pleasure.
\medskip

\section{Notation and background material}
In this section we fix notation and gather some background material
on Galois representations that is
important for the rest of this paper.

Let $X$ be a scheme. For $x\in X$ we denote by $k(x)$ the residue field at $x$.
If
$X$ is integral, then $R(X)$ stands for the function field of $X$, that is, for the
residue field at the generic point of $X$.
If $X$ happens to be a scheme of finite type over a base field $F$, then we often
write $F(x)$ instead of $k(x)$ and $F(X)$ instead of $R(X)$.

If $K$ is a field, then we denote by $K\sep$ (resp. $\widetilde{K}$) the
separable (resp. algebraic) closure of $K$ and by $G_K$ its absolute
Galois group. 
 A finitely generated field is by definition a field which is
finitely generated over its prime field.
For an abelian variety $A/K$ we let $\End_K(A)$ be the ring
of all $K$-endomorphisms of $A$. We denote by $\End(A):=\End_{\widetilde{K}}(A_{\widetilde{K}})$
the {\em absolute} endomorphism ring.

If $\Gamma$ is an object in an abelian category and $n\in\Zz$, then $n_\Gamma: \Gamma\to \Gamma$
is the morphism ``multiplication by $n$'' and $\Gamma[n]$ is the kernel of $n_\Gamma$.
Recall that there is an equivalence of categories between the category of finite \'etale
group schemes over $K$ and the category of finite (discrete) $G_K$-modules,
where we attach $\Gamma(K\sep)$ to a finite \'etale group scheme $\Gamma/K$.
For such a finite \'etale group scheme $\Gamma/K$ we sometimes write just
$\Gamma$ instead of $\Gamma(K\sep)$, at least in situations where we are sure that this
does not cause any confusion. 
For example,
if $A/K$ is an abelian variety and $n$ an integer coprime to $ \chara(K)$, then we often
write $A[n]$ rather than $A(K\sep)[n]$. Furthermore we put $A[n^\infty]:=\bigcup_{i\in\Nn} A[n^i]$.

If $M$ is a $G_K$-module (for example $M=\mu_n$ or $M=A[n]$ where
$A/K$ is an abelian variety), then we shall denote
the corresponding representation
of the Galois group $G_K$ by
$$\rho_M: G_K\to \Aut(M)$$
and define $\Mon_K(M):=\rho_M(G_K)$. We define $K(M):=K\sep^{\ker(\rho_M)}$ to be the fixed field in
$K\sep$ of the kernel of $\rho_M$. Then $K(M)/K$ is a Galois extension and
$G(K(M)/K)\cong \Mon_K(M)$.

If $R$ is a commutative ring with $1$ (usually $R=\Ff_\ell$ or $R=\Zz_\ell$) and $M$ is a finitely generated
free $R$-module equipped with a non-degenerate
alternating bilinear pairing $e:M\times M\to R'$ into
a free $R'$-module of rank $1$ (which is a multiplicatively written $R$-module
in our setting
below), then we denote by
$$\Sp(M, e)=\{f\in \Aut_R(M)\quad|\quad\forall x, y\in M: e(f(x), f(y))=e(x, y)\}$$
the corresponding symplectic group and by
$$\GSp(M, e)=\{f\in \Aut_R(M)\quad|\quad \exists\varepsilon\in R^\times:\forall x, y\in M:
e(f(x), f(y))=\varepsilon e(x, y)\}$$
the corresponding group of symplectic similitudes.

Let $n$ be an integer coprime to $ \chara(K)$ and $\ell$ be a prime different from
$ \chara(K)$. Let $A/K$ be an abelian variety. We denote by $A^\vee$ the dual abelian
variety and let $e_n: A[n]\times A^\vee[n]\to \mu_n$ and $e_{\ell^\infty}: T_\ell A\times
T_\ell A^\vee\to \Zz_\ell(1)$ be the corresponding Weil pairings. If $\lambda: A\to A^\vee$
is a polarization, then we deduce Weil pairings
$e_n^\lambda: A[n]\times A[n]\to \mu_n$ and $e_{\ell^\infty}^\lambda:
T_\ell A\times T_\ell A\to \Zz_\ell(1)$ in the obvious way.
If $\ell$ does not divide $\deg(\lambda)$ and if $n$ is coprime to $\deg(\lambda)$,
then $e_n^\lambda$ and $e_{\ell^\infty}^\lambda$ are non-degenerate, alternating,
$G_K$-equivariant pairings. Hence we have representations
$$\rho_{A[n]}: G_K\to \GSp(A[n], e_n^\lambda),$$
$$\rho_{T_\ell A}: G_K\to \GSp(T_\ell A, e_{\ell^\infty}^\lambda)$$
with images
$\Mon_K(A[n])\subset \GSp(A[n], e_n^\lambda)$ and
$\Mon_K(T_\ell A)\subset \GSp(T_\ell A, e_{\ell^\infty}^\lambda)$.
We shall say that an abelian variety $(A, \lambda)$ over a field $K$
has {\sl big monodromy}, if there is a constant $\ell_0>\max( \chara(K), \deg(\lambda))$ such
that $\Mon_K(A[\ell])\supset \Sp(A[\ell], e_\ell^\lambda)$ for every prime number
$\ell\ge \ell_0$.

Now let $S$ be a noetherian regular $1$-dimensional connected scheme  with function field $K=R(S)$ and $A/K$ an
abelian variety. Denote by $\AA\to S$ the N\'eron model (cf. \cite{blr}) of $A$. For $s\in S$
let $A_s:=\AA\times_S\Spec(k(s))$ be the corresponding fiber. Recall that we say that
$A$ has {\sl good reduction at $s$} provided $A_s$ is an abelian variety. In general,
we denote by $A_s^\circ$ the connected component of $A_s$. If $T$ is a maximal torus in $A_s^\circ$,
then $\dim(T)$ does not depend on the choice of $T$ \cite[IX.2.1]{SGA7} and we call $\dim(T)$ the
{\sl toric dimension} of the reduction $A_s$ of $A$ at $s$. Finally recall that one
says that $A$ has {\sl semi-stable reduction at $s$}, if $A_s^\circ$ is an extension
of an abelian variety by a torus.

We shall also need the following connections between the reduction type of $A$ and
properties of the Galois representations attached to $A$.
Let $s$ be a closed point of $S$. The valuation $v$ attached to $s$ admits an
extension to the separable closure $K\sep$;  we choose such an extension $\ol{v}$ and denote
by $D(\ol{v})$ the corresponding decomposition group. This is the
absolute Galois group of the quotient field $K_s=Q(\OO_{S, s}^h)$ of the henselization $\OO_{S, s}^h$
of the valuation ring $\OO_{S, s}$ of $v$. Hence the results mentioned in \cite[I.0.3]{SGA7}
for the henselian case carry over to give the following description of
$D(\ol{v})$:
If $I(\ol{v})$ is the kernel of the canonical map $D(\ol{v})\to G_{k(s)}$ defined by $\ol{v}$,
then $D(\ol{v})/I(\ol{v})\cong G_{k(s)}$. Let $p$
be the characteristic of the residue field $k(s)$ ($p$ is zero or a prime number). $I(\ol{v})$ has
a maximal pro-$p$ subgroup $P(\ol{v})$ ($P(\ol{v})=0$ if $p=0$) and
$$I(\ol{v})/P(\ol{v})\cong \plim_{n\notin p\Zz} \mu_n(k(s)\sep)\cong \prod_{\ell\neq p\prm} \Zz_\ell(1).$$
Hence the maximal pro-$\ell$-quotient $I_\ell(\ol{v})$ of $I(\ol{v})$ is isomorphic to $\Zz_\ell(1)$,
if $\ell\neq p$ is a prime.


\begin{prop} \label{semistablered} Let $\ell\neq p$ be a prime number. Assume that $A$ has semi-stable reduction at
$s$.
\begin{enumerate}
\item[a)] The image $\rho_{A[\ell]}(P(\ol{v}))=\{Id\}$ and $\rho_{A[\ell]}(I(\ol{v}))$
is a cyclic $\ell$-group.
\item[b)] Let $g$ be a generator of $\rho_{A[\ell]}(I(\ol{v}))$. Then $(g-Id)^2=0$.
\item[c)] Assume that $\ell$ does not divide the order of the component group
of $A_s$. The toric dimension of $A$ at $s$ is equal to $2\dim(A)-\dim_{\Ff_\ell}(\Eig(g, 1))$ if
$\Eig(g, 1)=\ke(g-Id)$ is the eigenspace of $g$ at $1$.
\end{enumerate}
\end{prop}

{\em Proof.} Part a) and b) are immediate consequences of \cite[IX.3.5.2.]{SGA7}.

Assume from now on that $\ell$ does not divide the order of the component group
of $A_s$. This assumption implies
$A_s^\circ [\ell]\cong A_s[\ell]$.

As we assumed $A$ to be semi-stable at $s$, there is an exact sequence
$$0\to T\to A_s^\circ \to B\to 0$$
where $T$ is a torus and $B$ is an abelian variety and $\dim(T)+\dim(B)=\dim(A_s)=\dim(A)$.
Now $\dim_{\Ff_\ell}(T[\ell])=\dim(T)$ and
$\dim_{\Ff_\ell}(B[\ell])=2\dim(B)=2\dim(A)-2\dim(T)$. Taking into account that we have an exact sequence
$$0\to T[\ell]\to A_s^\circ[\ell] \to B[\ell]\to 0$$
(note that $T(\widetilde{k})\cong (\widetilde{k}^\times)^{\dim(T)}$ is
divisible by $\ell$), we find the relation $\dim_{\Ff_\ell}(A_s[\ell])=\dim_{\Ff_\ell}(A_s^\circ[\ell])=2\dim(A)-\dim(T)$.
This implies c), because $A_s[\ell]=A[\ell]^{I(\ol{v})}$ (\cite[p. 495]{serretate1968}) and obviously
$A[\ell]^{I(\ol{v})}=\Eig(g, 1)$ .\hfill $\Box$

In general, if $V$ is a finite dimensional vector space over $\Ff_\ell$, and $g\in \End_{\Ff_\ell}(V)$, then
one defines $\dro(g)=\dim(V)-\dim(\Eig(g, 1))$. One calls $g$ a {\sl transvection}, if it is
unipotent of drop $1$. We shall say that an abelian variety $A$ over
a field $K$ is {\sl of Hall type}, provided $\End(A)=\Zz$ and there is
a discrete valuation $v$ on $K$ such that $A$ has semistable reduction of
toric dimension $1$ at $v$ (i.e. at the maximal ideal of the discrete valuation ring of $v$).
We have thus proved the following

\begin{prop} \label{extransv} If $A$ is an abelian variety of Hall type
over a finitely generated field $K$, then  there is a constant $\ell_0$ such that
$\Mon_K(A[\ell])$ contains a transvection for every prime number $\ell\ge \ell_0$.
\end{prop}

\section{Finiteness properties of division fields}
If $A$ is an abelian variety over a field $K$ (of arbitrary characteristic)
and $p= \chara(K)$,
then we denote by $A_{\neq p}$
the group of points in $A(K\sep)$ of order prime to $p$. Then
$$K(A_{\neq p})=\prod_{\ell\neq p\ prime} K(A[\ell^\infty])=\bigcup\limits_{n\notin p\Zz} K(A[n]).$$
If $p=0$, then $K(A_{\neq p})=K(A\tor)$. In this section we prove among other things: If
$K$ is finitely generated of positive characteristic, then $G(K(A_{\neq p})/K)$ is a
finitely generated profinite group.
We follow preprint \cite{fjp2010} as far as Lemmas \ref{fglemm1} and \ref{fgen} are concerned,
providing details of proofs for the reader's convenience.

In this section, a {\sl function field of $n$ variables}
over a field $F$ will be a finitely generated field extension $E/F$ of transcendence degree $n$. As usual
we call such a function field $E/F$ of $n$ variables {\sl separable} if it has a separating transcendency base.

\begin{lemm} Let $F$ be a separably closed field and $K/F$ a function field of
one variable. Assume that $K/F$ is separable. Put $p= \chara(F)$. Let $A/K$ be an abelian variety.
Then $G(K(A_{\neq p})/K)$ is
a finitely generated profinite group.\label{fglemm1}
\end{lemm}

{\em Proof.} There is a smooth projective curve $C/F$ with function field $K$. By Grothendieck's
Theorem \cite[IX.3.6]{SGA7} there is a finite separable extension $K'/K$ such that $A_{K'}$ has
semistable reduction at all points of the normalization $C'$ of $C$ in $K'$. We may assume that $K'/K$ is Galois.

Let $S'\subset C'$ be the finite set of closed points where $A_{K'}$ has bad reduction. Then for every $\ell\neq p$
the extension $K'(A[\ell^\infty])/K'$ is tamely ramified at all points of $C'$ (cf. \cite[IX.3.5.2.]{SGA7}) and
unramified outside $S'$ by the criterion of N\'eron-Ogg-Shafarevich \cite[Thm. 1]{serretate1968}. Hence
$K'(A_{\neq p})$ is contained in the maximal tamely ramified extension $K'_{S', tr}$ of $K'$ which is unramified
outside $S'$. The Galois group $G(K'_{S', tr}/K')$ is finitely generated by
\cite[Corollaire XIII.2.12]{SGA1}. Hence $G(K'(A_{\neq p})/K')$ is finitely generated as a quotient of
$G(K'_{S', tr}/K')$. Furthermore there is an exact sequence
$$1\to G(K'(A_{\neq p})/K')\to G(K(A_{\neq p})/K) \to G(K'/K)$$
and $G(K'/K)$ is finite. Hence $G(K(A_{\neq p})/K)$ is finitely generated as desired.\hfill $\Box$

\begin{lemm} Let $F$ be a field and $K/F$ a function field of one variable.
Assume that $K/F$ is separable.
Let $p= \chara(F)$. Let $A/K$ be an abelian variety. Let $F'$ be the
algebraic closure of $F$ in $K(A_{\neq p})$. Then $G(K(A_{\neq p})/F'K)$ is
a finitely generated profinite group.\label{fgen}
\end{lemm}

{\em Proof.}
$F\sep$ is $F'$-linearly disjoint from $K(A_{\neq p})$. Hence $F\sep K$ is
$F'K$-linearly disjoint from $K(A_{\neq p})$. This implies
$$G(K(A_{\neq p})/F'K)\cong G(F\sep K(A_{\neq p})/F\sep K),$$ and the latter
group is finitely generated by Lemma \ref{fglemm1} above.\hfill $\Box$

\begin{lemm} Let $(K, v)$ be a discrete valued field, $A/K$ an abelian variety with
good reduction at $v$, $n$ an integer coprime to the residue characteristic of $v$,
$L=K(A[n])$ and $w$ an extension of $v$ to $L$. Denote the residue field
of $v$ (resp. $w$) by $k(v)$ (resp. $k(w)$). Let $A_v/k(v)$ be the reduction
of $A$ at $v$. Then $k(w)=k(v)(A_v[n])$.\label{redlemm}
\end{lemm}

{\em Proof.} Let $R$ be the valuation ring of $v$ and $S=\Spec(R)$. Let $\AA\to S$ be an
abelian scheme with generic fibre $A$. Then $\AA[n]$ is a finite \'etale
group scheme over $S$. Let $T$ be the normalization
of $S$ in $L$. The restriction map $r: A[n](L)\cong \AA[n](T)\to A_v[n](k(w))$ is injective
\cite{serretate1968} and $|A[n](L)|=n^{2\dim(A)}$. Hence $r$ is an isomorphism and we may identify $A[n]$ with
$A_v[n]$. The fact that the whole $n$-torsion of $A_v$ is defined over $k(w)$ implies
that $k(v)(A_v[n])\subset k(w)$. We have to prove the other inclusion:
Let $D(w)$ be the decomposition group of the prime $w$ over $v$,
i.e. the stabilizer of $w$ under the
action of $G(L/K)$. Then $D(w)\to G(k(w)/k(v))$ is an isomorphism by
the criterion of N\'eron-Ogg-Shafarevich. As $D(w)\to Aut(A[n])$ is
injective, it follows that $G(k(w)/k(v))\to Aut(A_v[n])$ is injective
as well. This implies that $k(v)(A_v[n])= k(w)$.\hfill $\Box$

\begin{defi}
We shall say in the sequel that a field $K$ has property $\FF$, if $G(K'(A_{\neq p})/K')$ is
a {\em finitely generated} profinite group for every finite separable extension $K'/K$ and
every abelian variety $A/K'$.
\end{defi}

\begin{lemm} Let $F$ be a field that has property $\FF$. Let $p= \chara(F)$.
Let $K/F$ be a
function field of one variable. Assume that $K/F$ is separable. Then
$K$ has property $\FF$.\label{finallemm}
\end{lemm}

{\em Proof.} We have to show that $G(K'(A_{\neq p})/K')$ is finitely
generated for every finite separable extension $K'/K$ and every abelian
variety $A/K'$. But if $K'/K$ is a finite separable extension, then
$K'/F$ is a separable function field of one variable again. Hence it is
enough to prove that $G(K(A_{\neq p})/K)$ is finitely generated
for every abelian variety $A/K$.

Let $A/K$ be an abelian variety. Let $F_0$ be the algebraic closure
of $F$ in $K$. Then $K/F_0$ is a regular extension. Let $C/F_0$ be a smooth
curve with function field $K$ and such that $A$ has good reduction at
all points of $C$. There is a finite Galois extension $F_1/F_0$
such that $C(F_1)\neq \emptyset$. If we put $K_1:=F_1 K$, then
$K_1/F_1$ is regular. Furthermore
there is an exact sequence
$$1\to G(K_1(A_{\neq p})/K_1)\to G(K(A_{\neq p})/K)\to G(K_1/K)$$
and $G(K_1/K)$ is finite. If we prove that $G(K_1(A_{\neq p})/K_1)$
is finitely generated, then it follows that
$G(K(A_{\neq p})/K)$ is finitely generated as well.
Hence we may assume that $K_1=K$, i.e. that $K/F$ is regular and that
$C(F)\neq \emptyset$.

Choose a point $c\in C(F)$ and denote by $A_c/F$ the (good) reduction
of $A$ at $c$. As in Lemma \ref{fgen} denote by $F'$ the algebraic
closure of $F$ in $K(A_{\neq p})$.

{\bf Claim.} $F'\subset F(A_{c, \neq p})$.

Let $x\in F'$. Then $x$ is algebraic over $F$ and $x\in K(A[n])$ for
some $n$ which is coprime to $p$. If $F_n$ denotes the algebraic
closure of $F$ in $K(A[n])$, then $x\in F_n$. Let $w$ be the extension to $K(A[n])$ of
the valuation attached to $c$. Then $k(w)=F(A_c[n])$ by Lemma \ref{redlemm}.
Obviously $F_n\subset k(w)$. Hence $x\in F(A_c[n])\subset F(A_{c, \neq p})$.
This finishes the proof of the Claim.

The profinite group $G(F(A_{c, \neq p})/F)$ is finitely generated, because
$F$ has property $\FF$ by assumption.
Hence its quotient $G(F'/F)$ is finitely generated as well. Note that
$G(F'K/K)=G(F'/F)$. On the
other hand $G(K(A_{\neq p})/F'K)$ is finitely generated by
Lemma \ref{fgen}. From the exact sequence
$$1\to G(K(A_{\neq p})/F'K)\to G(K(A_{\neq p})/K)\to G(F'K/K)\to 1$$
we see that $G(K(A_{\neq p})/K)$ is finitely generated as desired.\hfill $\Box$

\begin{prop} Let $F$ be a perfect field which has property $\FF$. Then
every finitely generated extension $K$ of $F$ has property $\FF$.
\end{prop}

{\em Proof.} We prove this by induction on $\trdeg(K/F)$. If $\trdeg(K/F)=0$ there
is nothing to prove. Assume $\trdeg(K/F)=d\ge 1$. We may assume that every
finitely generated extension $F'$ of $F$ with $\trdeg(F'/F)<d$ has property $\FF$.

Choose a separating transcendency base $(x_1,\cdots, x_d)$ for $K/F$. Put
$F':=K(x_1,\cdots, x_{d-1})$. Then $F'$ has property $\FF$ by the induction
hypothesis. Furthermore $K/F'$ is a function field of one variable and
$K/F'$ is separable. Hence Lemma \ref{finallemm} implies that $K$ has property $\FF$.\hfill $\Box$

\begin{coro} Let $K$ be a finitely generated field {\em of positive
characteristic} or $K$ be a function field over an algebraically closed field of arbitrary characteristic.
Then $K$ has property $\FF$. In particular $G(K(A_{\neq p})/K)$
is finitely generated for every abelian variety $A/K$. \label{fgl}
\end{coro}

{\em Proof.} A finite field $\Ff$ is perfect. It has property $\FF$, because its
absolute Galois group is procyclic. An algebraically closed field is perfect and has property $\FF$,
because its absolute Galois group is the trivial group. Now in both cases
$K$ is a function field over a perfect field which has property $\FF$.\hfill $\Box$

\begin{rema} A finitely generated field $K$ of characteristic zero does {\em not}
have property $\FF$. In fact, if $A/K$ is principally polarized
abelian variety, then by the existence of the Weil pairing $K(A\tor)\supset K(\mu_\infty),$ and
plainly $G(K(\mu_\infty)/K)$ is {\em not} finitely generated, when $K$ is a finitely generated extension
of $\Qq$.\label{fglrema}
\end{rema}

\section{Monodromy Computations}
Let $K$ be a field and $A/K$ an abelian variety. We begin with the question whether $A[\ell]$ is a simple
$G_K$-module for sufficiently large $\ell$. In the cases we need to consider,
this question has an affirmative answer due to the following classical fact (cf. \cite[p. 118, p. 204]{fwbook},
 \cite{zarhin1977}, \cite{zarhin1985},\cite{mb1985}).

\begin{thm} {\bf (Faltings, Zarhin)} \label{fz}
Let $K$ be a finitely generated field and $A/K$ an abelian variety.
Then there is a constant $\ell_0> \chara(K)$ such that the $\Ff_\ell[G_K]$-module $A[\ell]$ is semisimple and the canonical map
$\End_K(A)\otimes \Ff_\ell\to \End_{\Ff_\ell}(A[\ell])$ is injective with image
$\End_{\Ff_\ell[G_K]}(A[\ell])$ for all primes $\ell\ge \ell_0$.
\end{thm}

\begin{prop} \label{irr} Let $A$ be an abelian variety over a finitely generated field $K$.
Assume that $\End_K(A)=\Zz$. Then there is a constant $\ell_0$ such that
$A[\ell]$ is a simple $\Ff_\ell[G_K]$-module for all primes $\ell\ge \ell_0$.
\end{prop}

{\em Proof.} By Theorem \ref{fz} there is a constant $\ell_0$ such that $A[\ell]$ is a semisimple $\Ff_\ell[G_K]$-module
with $\End_{\Ff_\ell[G_K]}(A[\ell])=\Ff_\ell Id$
for every prime $\ell\ge \ell_0$. This is only possible if $A[\ell]$ is a simple $\Ff_\ell[G_K]$-module
for all primes $\ell\ge \ell_0$.\hfill $\Box$

We need some notation in order to explain a theorem of Raynaud that will be of importance later. 
Let $E/\Ff_p$ be a finite field extension with $|E|=p^d$ and $F/\Ff_p$ an algebraic extension.
Denote by $\Emb(E, \widetilde{F})$ the set of all embeddings $E\to \widetilde{F}$. If $i\in \Emb(E, \widetilde{F})$ is one such
embedding, then $\Emb(E, \widetilde{F})=\{i^{p^a}: a\in\{0,\cdots, d-1\}\}$. Furthermore the restriction
$i|E^\times$ lies in $\Hom(E^\times, \widetilde{F}^\times)$. For
every character $\chi\in \Hom(E^\times, F^\times)$ there is a unique $m\in \{0, \cdots, p^d-2\}$ such
that $\chi=(i|E^\times)^m$. Expanding $m$ $p$-adically, we see that there is a unique function
$e: \Emb(E, \widetilde{F})\to \{0, \cdots, p-1\}$ such that  
$$\chi=\prod_{j\in \Emb(E, \widetilde{F})} (j|E^\times)^{e(j)},$$
and such that $e(j)<p-1$ for some $j\in\Emb(E,\widetilde{F})$.
We define $\amp(\chi):=\max(e(j): j\in \Emb(E, \widetilde{F}))$ to be the {\em amplitude of the character $\chi$}. 
Let $\rho: E^\times\to\Aut_{\Ff_p}(V)$ be a representation of $E^\times$ on a finite dimensional
$\Ff_p$-vector space $V$. If $V$ is a {\em simple} $\Ff_p[E^\times]$-module, then there is 
a finite field $F_V$ with $|F_V|=|V|$ and a structure of $1$-dimensional $F_V$-vector space
on $V$ such that $\rho$ factors through a character $\chi_\rho: E^\times\to F_V^\times$. We then
define $\amp(V):=\amp(\rho):=\amp(\chi_\rho)$. In general $V$ is a semisimple $\Ff_p[E^\times]$-module
by Maschke's theorem, and we can write $V=V_1\oplus\cdots\oplus V_t$ 
as a direct sum of simple $\Ff_p[E^\times]$-modules and define $\amp(V):=\amp(\rho):=\max(\amp(V_i): i=1, \cdots, t)$
to be the {\em amplitude of the representation $\rho$}.  
With this terminology in mind, we can state Raynaud's theorem in the following way.

\begin{thm} {\bf (Raynaud \cite{raynaud1974}, \cite[p. 277]{serre1972})} Let $A$ be an abelian variety over a number field $K$. Let $v$ be a 
place of $K$ with residue characteristic $p$. Let $e$ be the ramification index of $v|\Qq$.
Let $w$ be an extension of $v$ to $K(A[p])$. Let 
$I$ be the inertia group of $w|v$ and $P$ the $p$-Sylow subgroup of $I$. Let $C\subset I$ be 
a subgroup that maps isomorphically onto $I/P$. Then there is a finite extension $E/\Ff_p$ and
a surjective homomorphism $E^\times\to C$ such that the resulting representation
$$\rho: E^\times\to C\to \Aut_{\Ff_p}(A[p])$$
has amplitude $\amp(\rho)\le e$.\label{raynaud}
\end{thm}

The technical heart of our monodromy computations is the following
group theoretical result, which can be extracted from
the work of C. Hall \cite{hall2007}, \cite{hall2008}. 

\begin{thm} \label{stronghall} Let $\ell>2$ be a prime, let $(V, e_V)$ be a finite-dimensional symplectic space over
$\Ff_\ell$ and $M$ a subgroup of $\Gamma:=\GSp(V, e_V)$. Assume that $M$ contains
a transvection and that $V$ is a simple $\Ff_\ell[M]$-module. Denote by
$R$ the subgroup of $M$ generated by the transvections in $M$.
\begin{enumerate}
\item[a)] Then there is a non-zero symplectic subspace $W\subset V$, which is
a simple $\Ff_\ell[R]$-module, such that the following properties hold true:
\begin{enumerate}
\item[i)] Let $H=\Stab_M(W)$. There is a orthogonal direct
sum decomposition $V=\bigoplus_{g\in M/H} gW$. In particular $|M/H|\le \dim(V)$.
\item[ii)] $R\cong \prod_{g\in M/H} \Sp(W)$ and $N_\Gamma(R)\cong \prod_{g\in M/H} \GSp(W)\rtimes
\Sym(M/H)$.
\item[iii)] $R\subset M\subset N_\Gamma(R)$.
\end{enumerate}
Denote by $\varphi: N_\Gamma(R)\to \Sym(M/H)$ the projection.
\item[b)] Let $e\in \Nn$.
Let $E/\Ff_\ell$ be a finite extension and $\rho: E^\times\to M\subset \GSp(V,e_V)$ a homomorphism such
that the corresponding representation of $E^\times$ on $V$ has amplitude $\amp(\rho)\le e$. 
If $\ell>\dim(V)e+1$, then $\varphi(\rho(E^\times))=\{1\}.$
\end{enumerate}
\end{thm}

Hall's proof in \cite{hall2007}, \cite{hall2008}
adresses a slightly 
less general situation. We will present a self-contained proof of Theorem
\ref{stronghall} in the Appendix.

\begin{rema} \label{hallrema} Assume that in the situation of Theorem \ref{stronghall}
the module $V$ is a simple $\Ff_\ell[\ke(\varphi)\cap M]$-module. Then $V$ is in
particular a simple
$\Ff_\ell[\ke(\varphi)]$-module and $\ke(\varphi)=\prod_{g\in M/H} \GSp(W)$. This 
is only possible if $M=H$, $V=W$ and $R=\Sp(V, e)\subset M$.
\end{rema}

We now state the main result of this section.

\begin{thm} Let $K$ be a finitely generated field. Let
$(A, \lambda)$ be a polarized abelian variety over $K$ of Hall type.
Then $(A, \lambda)$ has big monodromy. \label{main1}
\end{thm}

The case where $K$ is a global field is due to Hall (cf. \cite{hall2008}) and we
follow his line of proof to some extent, but we need a lot of additional arguments
in order to make things work in the more general situation. The proof will occupy almost the rest of this section.

There is a constant $\ell_0>\max(\deg(\lambda), \chara(K))$ such that the following holds true for
all primes $\ell\ge \ell_0$:

\begin{enumerate}
\item The subgroup $\Mon_K(A[\ell])$ of $\GSp(A[\ell], e_\ell^\lambda)$
contains a transvection. Denote by $R_\ell$ the subgroup of
$\Mon_K(A[\ell])$ generated by the transvections in $\Mon_K(A[\ell])$ (cf.
Proposition \ref{extransv}).
\item $A[\ell]$ is a simple $\Ff_\ell[G_K]$-module (cf. Proposition \ref{irr}).
\end{enumerate}

Now Hall's group theory result (cf. Theorem \ref{stronghall})
gives - for every prime $\ell\ge \ell_0$ - a non-zero
symplectic subspace $W_\ell\subset A[\ell]$, which is simple as a $\Ff_\ell[R_\ell]$-module such that
the properties i), ii) and iii) of Theorem \ref{stronghall} are satisfied. Let $H_\ell$
be the stabilizer of $W_\ell$ under the action of $\Mon_K(A[\ell])$. Define $M_\ell:=
\Mon_K(A[\ell])$ and $\Gamma_\ell:=\GSp(A[\ell], e_\ell^\lambda)$.
Then
$$\prod_{M_\ell/H_\ell} \Sp(W_\ell, e_\ell^\lambda)\cong R_\ell\subset M_\ell\subset N_{\Gamma_\ell}(R_\ell)=
\prod_{M_\ell/H_\ell} \Sp(W_\ell, e_\ell^\lambda)\rtimes \Sym(M_\ell/H_\ell),$$
and we denote by $\varphi_\ell: N_{\Gamma_\ell}(R_\ell)\to \Sym(M_\ell/H_\ell)$ the projection. We
have the following property (cf. Remark \ref{hallrema}):

{\em If $A[\ell]$ is a simple $\Ff_\ell[\ke(\varphi_\ell)\cap M_\ell]$-module for some
prime $\ell\ge \ell_0$, then
$M_\ell=H_\ell$, $W_\ell=A[\ell]$ and $M_\ell\supset \Sp(A[\ell], e_\ell^\lambda)$ for this prime $\ell$.}

We denote by $N_\ell$ the fixed field inside $K\sep$ of
the preimage $\rho_{A[\ell]}^{-1}(M_\ell\cap \ke(\varphi_\ell))$, where
$\rho_{A[\ell]}: G_K\to \Gamma_\ell$ is the mod-$\ell$ representation attached to $A$.
Then $N_\ell$ is an intermediate field of $K(A[\ell])/K$ which is Galois over $K$,
and $G(N_\ell/K)$ is isomorphic to the subgroup
$\varphi_\ell(M_\ell)$ of $\Sym(M_\ell/H_\ell)$. In particular $[N_\ell:K]\le (2\dim(A))!$ is
bounded independently of $\ell$. If we denote by $N:=\prod_{\ell\ge \ell_0\prm} N_\ell$ the
corresponding composite field, then $G_N=\bigcap_{\ell\ge \ell_0\prm} G_{N_\ell}$.
Hence the following property holds true.

{\em If $A[\ell]$ is simple as a $\Ff_\ell[G_N]$-module for some prime $\ell\ge \ell_0$, then
$M_\ell\supset \Sp(A[\ell], e_\ell^\lambda)$ for this prime $\ell$.}\hfill $(*)$

{\em Proof of Theorem \ref{main1} in the special case $\chara(K)>0$.}
If $\chara(K)>0$, then the Galois group $G(K(A_{\neq p})/K)$ ($p:=\chara(K)$) is finitely generated, because $K$ then has property $\FF$ by Corollary
\ref{fgl}. Furthermore $N_\ell$ is an intermediate field of
$K(A_{\neq p})/K$ which is Galois over $K$ and with $[N_\ell:K]$
bounded independently of $\ell$.
Hence $N/K$ must be finite. In particular $N$ is finitely generated.
A second application of the result of Faltings and Zarhin (cf.
Proposition \ref{irr}) yields a constant $\ell_1\ge \ell_0$ such that
$A[\ell]$ is a simple $\Ff_\ell[G_N]$-module for all primes $\ell\ge \ell_0$. Hence
$A$ has big monodromy by $(*)$.\hfill $\Box$

To finish the proof of Theorem \ref{main1} we assume for the rest of the proof 
that $\chara(K)=0$. We shall prove that $N/K$ is finite also in that case, but
the proof of this fact is more complicated, because now $K$ is {\em not}
$\FF$-finite (cf. Remark \ref{fglrema}). We briefly sketch the main steps in
the proof, before we go into the details: The first and hardest step is to show
that the algebraic closure $L$ of $\Qq$ in $N$ is a {\em finite} extension
of $\Qq$. In order to achieve this we will construct a
finite extension $L'/\Qq$ such that some $L'$-rational
``place'' of $KL'$ splits up completely into $L'$-rational ``places'' of
$N_\ell L'$ for every sufficiently large
prime $\ell$. We use this to show that $G(NL/KL)\cong
G(NL\sep/KL\sep)$ and the fact that the latter group can be proved to be finite, because
$KL\sep$ is $\FF$-finite (unlike $K$ itself). This suffices to prove that $N/K$ is
finite.  Once we know this, we shall proceed as in the positive characteristic
case above.

We now go into the details. Let $F$ be the algebraic closure
of $\Qq$ in $K$. Then $F$ is a number field.
Let $S$ be a smooth affine $F$-variety with function field
$K$ such that $A$ extends to an abelian scheme $\AA$ over $S$ with generic fibre $A$
(i.e. such that $A$ has good reduction along $S$). Let $S_\ell$ be the normalization
of $S$ in $N_\ell$ and let $S_\ell'$ be the normalization of $S_\ell$ in $K(A[\ell])$.
Then $S_\ell'\to S_\ell\to S$ are finite \'etale covers. (Note that $\chara(F(s))=0$ for every point $s\in S$.)
In particular $S_\ell'$ and $S_\ell$ are smooth $F$-schemes. (Compare the diagram below.)

Fix a geometric point $P\in S(F\sep)$ and denote by $A_P:=\AA\times_S \Spec(F(P))$
the corresponding special fibre of $\AA$. Then $A_P$ is an abelian variety over
the number field $F(P)$.
Fix for every $\ell\ge \ell_0$ a geometric point $Q_\ell\in S_\ell(F\sep)$ over $P$ and a geometric
point
$Q_\ell'\in S_\ell'(F\sep)$ over $Q_\ell$. Then $F(Q_\ell')/F(Q_\ell)$ and $F(Q_\ell)/F(P)$ are
finite extensions of number fields. Note that $F(Q_\ell')=F(P)(A_P[\ell])$ by
Lemma \ref{redlemm}.
Denote by $\OO$ (resp. $\OO_\ell$, resp. $\OO_\ell'$) the
integral closure of $\Zz$ in $F(P)$ (resp. in $F(Q_\ell)$, resp. in $F(Q_\ell')$). For
every prime $\ell\ge \ell_0$ we have the following diagram on the level of schemes

$$\begin{xy}
  \xymatrix{
&\Spec(K(A[\ell]))\ar[r]\ar[d]       & \Spec(N_\ell)\ar[r]\ar[d]         & \Spec(K)\ar[d]\\
&   S_\ell'\ar[r]           &  S_\ell\ar[r]             & S \\
\Spec(F(P)(A_P[\ell]))\ar@{=}[r]&\Spec(F(Q_\ell'))\ar[r]\ar[d]\ar[u]&\Spec(F(Q_\ell))\ar[r]\ar[d]\ar[u] & \Spec(F(P))\ar[d]\ar[u] \\
&\Spec(\OO_\ell')\ar[r]          &\Spec(\OO_\ell)\ar[r]^{f_\ell}     &  \Spec(\OO)
}
\end{xy}$$

We now study the ramification of prime ideals $\mf m\in \Spec(\OO)$
in the extension $F(Q_\ell)/F(P)$.
Let $\Pp\bad$ be the (finite) set of primes $\mf p\in
\Spec(\OO)$ where $A_P/F(P)$ has bad reduction.

\begin{lemm} \label{halllemm} There is a constant $\ell_2\ge \ell_0$ with the following
property: For every prime number $\ell\ge \ell_2$ the map $f_\ell: \Spec(\OO_\ell)\to
\Spec(\OO)$ is \'etale at every point $\mf m\in \Spec(\OO)$ outside of $\Pp\bad$.
\end{lemm}

{\em Proof.}
Let $\ell_2:=\max(\ell_0, (2\dim(A))![F(P):\Qq]+2)$.

Now let $\ell\ge \ell_2$ be a prime number. Let $\mf m\in \Spec(\OO)$ be
an arbitrary prime ideal with
$\mf m\notin \Pp\bad$. We have to show that $\mf m$ is unramified
in $F(Q_\ell)$. Let $p=\chara(\OO/\mf m)$ be the residue characteristic of $\mf m$.

If $p\neq \ell$, then $\mf m$ is unramified even in $F(Q_\ell')=F(P)(A_P[\ell])$.

We can hence assume that $\fbox{$p=\ell$}$. Let $\mf m_\ell\in \Spec(\OO_\ell)$ be
a point over $\mf m$ and $\mf m_\ell'\in \Spec(\OO_\ell')$ a point over $\mf m_\ell$.
Let $D(\mf m_\ell')$ (resp. $D(\mf m_\ell)$)
be the decomposition group of $\mf m_\ell'/F(P)$ (resp. of $\mf m_\ell/F(P)$)
and $I(\mf m_\ell')$ (resp. $I(\mf m_\ell)$) the
corresponding inertia group. Let $P(\mf m_\ell')$ (resp. $P(\mf m_\ell)$) be the (unique)
$p$-Sylow subgroup of $I(\mf m_\ell')$ (resp. $I(\mf m_\ell)$).

We have the following commutative diagram on the level of groups:

$$\begin{xy}
  \xymatrix{
 \prod_{M_\ell/H_\ell} \GSp(W_\ell)\ar@^{(->>}[r] & N_{\Gamma_\ell}(M_\ell)\ar@{->>}[r] & Sym(M_\ell/H_\ell)\\
  M_\ell\cap ker(\varphi_\ell)\ar@^{(->}[r]\ar@^{(->}[u]  & M_\ell\ar@{->>}[r]\ar@^{(->}[u] & \varphi_\ell(M_\ell)\ar@^{(->}[u]\\
G(K(A[\ell])/N_\ell)\ar@^{(->}[r]\ar@{=}[u]        &
G(K(A[\ell])/K)\ar@{->>}[r]\ar@{=}[u] & G(N_\ell/K)\ar@{=}[u] \\
G(F(Q_\ell')/F(Q_\ell))\ar@^{(->}[r]\ar@^{(->}[u] &
G(F(Q_\ell')/F(P))\ar@{->>}[r]\ar@^{(->}[u] &
G(F(Q_\ell)/F(P))\ar@^{(->}[u]\\
 & D(\mf m_\ell')\ar@{->>}[r]\ar@^{(->}[u] & D(\mf m_\ell)\ar@^{(->}[u]\\
 & I(\mf m_\ell')\ar@{->>}[r]\ar@^{(->}[u] & I(\mf m_\ell)\ar@^{(->}[u]\\
 & P(\mf m_\ell')\ar@{->>}[r]\ar@^{(->}[u] & P(\mf m_\ell)\ar@^{(->}[u]\\
}
\end{xy}$$

We have to prove that the image of $I(\mf m_\ell')$ in
 $\Sym(M_\ell/H_\ell)$ by the maps in the diagram is $\{1\}$.
Now $p=\ell> (2\dim(A))!$ due to our choice of $\ell_2$ and
$|\Sym(M_\ell/H_\ell)|\le (2\dim(A))!$, hence $P(\mf m_\ell')$ maps to $\{1\}$ in
$\Sym(M_\ell/H_\ell)$. In particular, $P(\mf m_\ell)=\{1\}$.
Consider the tame ramification group
$I_t=I(\mf m_\ell')/P(\mf m_\ell')$. It is a cyclic group of order
prime to $p$. Choose a subgroup $C\subset I(\mf m_\ell')$ that maps
isomorphically onto $I_t$ under the projection. It is enough to
show that $C$ maps to $\{1\}$ in $\Sym(M_\ell/H_\ell)$.

By Raynaud's theorem (cf. Theorem \ref{raynaud})
there is a finite extension $E/\Ff_p$ and an epimorphism $E^\times\to C$ such that the
resulting representation
$$E^\times\to C\to \Aut(A_P[\ell])= \Aut(A[\ell])$$
has amplitude $\le e$, where $e$ is the ramification index of $\mf m$ over $\Qq$. Clearly
$e\le [F(P):\Qq]$.
By part b) of Theorem \ref{stronghall}, the image of $E^\times$ in $\Sym(M_\ell/H_\ell)$ is $\{1\}$.
Hence the image of $C$ in  $\Sym(M_\ell/H_\ell)$ is $\{1\}$ as desired. $\hfill \Box$

\begin{lemm} \label{lemm2} Let $L$ be the algebraic closure of $F$ in $N$.
Then $L/F$ is a finite extension.
\end{lemm}

{\em Proof.} Let $L':=\prod_{\ell\ge \ell_0\prm} F(Q_\ell)$. For every
prime $\ell\ge \ell_2$ the Galois extension of number fields $F(Q_\ell)/F(P)$ is unramfied outside
$\Pp\bad$ by Lemma \ref{halllemm}. Furthermore
$[F(Q_\ell):F(P)]\le (2\dim(A))!$ for every prime $\ell\ge \ell_2$. The Theorem of Hermite-Minkowski
(cf. \cite{langANT}, p. 122) implies that $\prod_{\ell\ge \ell_2\prm} F(Q_\ell)$ is a {\em finite} extension of $F(P)$. This 
in turn implies that $L'/F$ is a finite extension. 
It is thus enough to show that $L\subset L'$.

Recall that $K=F(S)$ is the function field of the $F$-variety $S$ and $S_\ell$ is the normalization
of $S$ in the finite Galois extension $N_\ell/K$. Denote by $\hat{S}$ the normalization of $S$ in
$N$ and by $h_\ell: \hat{S}\to S_\ell$ the canonical projection. The
canonical morphism $\hat{S}\to S$ is surjective, hence there is a point $\hat{P}\in \hat{S}(F\sep)$
over $P$. The point $h_\ell(\hat{P})\in S_\ell(F\sep)$ lies over $P$. Hence $h_\ell(\hat{P})$ is
conjugate to $Q_\ell$ under the action of $G(N_\ell/K)$. This implies that $F(h_\ell(\hat{P}))=F(Q_\ell)$.
For every $\ell\ge \ell_0$ there is a diagram
$$\begin{xy}
  \xymatrix{
  \Spec(N)\ar[r]                 & \Spec(N_\ell)\ar[r]         & \Spec(K)\\ 
   \hat{S} \ar[r]^{h_\ell}\ar[d]\ar[u]          &  S_\ell\ar[r]\ar[d]\ar[u]         & S\ar[d]\ar[u] \\
   \Spec(F(\hat{P}))\ar[r]       &\Spec(F(Q_\ell))\ar[r]       & \Spec(F(P)) \\
}
\end{xy}$$
where the morphisms $S_\ell\to S$ are \'etale covers and $N=\prod_{\ell\ge\ell_0} N_\ell$. 
It follows that $F(\hat{P})=\prod_{\ell\ge\ell_0} F(Q_\ell)=L'$. On the other hand $L$ is the
algebraic closure of $F$ in $N$, hence $\hat{S}$ is a scheme over $L$. This implies that $L$ is
a subfield of $F(\hat{P})$. Hence in fact $L\subset L'$ as desired.\hfill $\Box$

{\em End of the proof of Theorem \ref{main1} in the case $\chara(K)=0$.}
We have an isomorphism $G(NL\sep/KL\sep)\cong G(N/KL)$, because
$N/L$ and $KL/L$ are regular extensions. The field $KL\sep$ is
$\FF$-finite by Corollary \ref{fgl}. Hence the profinite group $G(KL\sep(A\tor)/KL\sep)$
is finitely generated. As $NL\sep\subset KL\sep(A\tor)$, 
$G(NL\sep/KL\sep)$ must be finitely generated as well. 
Furthermore
$NL\sep =\prod_{\ell\ge \ell_0} N_\ell L\sep$ where $[N_\ell L\sep:KL\sep]$ is bounded independently from $\ell$.
Hence $G(NL\sep/KL\sep)$ is finite and this implies that $N/KL$ is a  finite extension.
On the other hand it follows from Lemma \ref{lemm2} that $KL/K$ is finite. Hence $N/K$ is 
a {\em finite} extension. Consequently 
$N$ is finitely generated, because $K$ is finitely generated. 
Proposition \ref{irr} yields a constant $\ell_3>\ell_0$ such
that $A[\ell]$ is a simple $\Ff_\ell(G_N)$-module for every prime $\ell\ge \ell_3$. Hence
$A/K$ has big monodromy by $(*)$, as desired.\hfill $\Box$

\section{Applications}

In this section we apply our methods to prove a generalization of a result of Ellenberg, Elsholz, Hall and Kowalski
on endomorphism rings and simplicity of fibres in certain families of abelian varieties (cf. \cite[Theorem 8]{ehk2009}).

\begin{prop} \label{monpropnew1}Let $K$ be a field and $(A, \lambda)$ a polarized
abelian variety over $K$ with big monodromy. Let $L/K$ be
a finite extension. Then the following properties hold.
\begin{enumerate}
\item[a)] There is a constant $\ell_0\ge \max(\chara(K), \deg(\lambda))$ such that $\Mon_L(A[\ell])\supset \Sp(A[\ell], e_\ell^\lambda)$
for every prime number $\ell\ge \ell_0$.
\item[b)] $A$ is geometrically simple.
\end{enumerate}
\end{prop}

{\em Proof.}
Part a). Let $E_0$ be the maximal separable extension of $K$ in $L$ and $E/K$ a finite Galois
extension containing $E_0$. By our assumption there is a
constant $\ell_0>\max(\deg(\lambda), \chara(K), 5)$ such that $\Mon_{K}(A[\ell])\supset\Sp(A[\ell], e_\ell^\lambda)$
for every prime $\ell\ge \ell_0$. For $\ell\ge \ell_0$ let $K_\ell$ be the fixed field of $\Sp(A[\ell], e_\ell^\lambda)$
in $K(A[\ell])/K$. Then $\Mon_{K_\ell}(A[\ell])=\Sp(A[\ell], e_\ell^\lambda)$ and
$\Mon_{EK_\ell}(A[\ell])$ is a normal
subgroup of $\Mon_{K_\ell}(A[\ell])$ 
of index $\le [E:K]$. Put $\ell_1:=\max(\ell_0, [E:K]+1)$. Then
$$|\Mon_{EK_\ell}(A[\ell])|\ge \frac{1}{[E:K]}|\Sp(A[\ell], e_\ell^\lambda)|>2$$
for all primes $\ell\ge \ell_1$.
On the other hand the only normal subgroups of $\Sp(A[\ell], e_\ell^\lambda)$ are $\{\pm 1\}$ and the
trivial group (cf. \cite[p. 53]{serretovigneras}). Hence $$\Mon_{E_0}(A[\ell])\supset\Mon_E(A[\ell])\supset
\Mon_{EK_\ell}(A[\ell])=\Sp(A[\ell], e_\ell^\lambda)$$
for all primes $\ell\ge \ell_1$. As $L/E_0$ is purely inseparable, we find
$$\Mon_L(A_L[\ell])=\Mon_{E_0}(A[\ell])\supset\Sp(A[\ell], e_\ell^\lambda)$$
for all primes $\ell\ge \ell_1$ as desired.  

Part b). Let $A_1, A_2/\widetilde{K}$ be abelian varieties and $f: A_{\widetilde{K}}\to A_1\times A_2$
an isogeny. Then $A_1, A_2$ and $f$ are defined over some finite extension $L/K$. Hence there
is an $\Ff_\ell[G_L]$-module isomorphism $A[\ell]\cong A_1[\ell]\times A_2[\ell]$ for every prime
$\ell>\deg(f)$. By part a) $\Mon_L(A[\ell])\supset \Sp(A[\ell]), e_\ell^\lambda)$ for all sufficiently large primes $\ell$.
Hence $A[\ell]$ is a simple $\Ff_\ell[\Mon_L(A[\ell])]$-module and in particular a simple $\Ff_\ell(G_L)$-module 
for all sufficiently large primes $\ell$.
This is only possible if $A_1=0$ or $A_2=0$.\hfill $\Box$

Let $F$ be a finitely generated field and
$K/F$ a finitely generated transcendental field extension and $A/K$ an
abelian variety. We say that $A/K$ is {\em weakly isotrivial with respect to $F$}, if there is an
abelian variety $B/\widetilde{F}$ and
a $\widetilde{K}$-isogeny $B_{\widetilde{K}}\to A_{\widetilde{K}}$.

\begin{prop} \label{geo} Let $F$ be a finitely generated field,
$K/F$ a finitely generated separable transcendental field extension and $(A, \lambda)$
a polarized abelian variety over $K$. Assume that $A/K$ has big monodromy
and that $A/K$ is {\em not} weakly isotrivial with respect to $F$.
Define $K':=F\sep K$. Then there is a constant $\ell_0\ge \max(\chara(K), \deg(\lambda))$ such that
$\Mon_{K'}(A[\ell])=\Sp(A[\ell], e_\ell^\lambda)$ for every prime number $\ell\ge \ell_0$.
\end{prop}

{\em Proof.} Let $\ell_0\ge \max(\deg(\lambda), \chara(K), 5)$ be a constant
such that $\Mon_K(A[\ell])\supset \Sp(A[\ell], e_\ell^\lambda)$ for every
prime $\ell\ge \ell_0$. Let $\ell\ge \ell_0$ be a prime number. Then we have
$$\Mon_{K'}(A[\ell])\subset \Sp(A[\ell], e_\ell^\lambda)\subset \Mon_K(A[\ell]),$$
because $K'$ contains $\mu_\ell$.
Furthermore $\Mon_{K'}(A[\ell])$ a normal subgroup of $\Mon_{K}(A[\ell])$, because
$K'/K$ is Galois.
It follows that $\Mon_{K'}(A[\ell])$ is {\em normal} in $\Sp(A[\ell], e_\ell^\lambda)$.

The only proper normal subgroups in
$\Sp(A[\ell], e_\ell^\lambda)$ are $\{1\}$ and $\{\pm 1\}$ (cf. \cite[p. 53]{serretovigneras}), because $\ell\ge 5$.
Hence either $\Mon_{K'}(A[\ell])=\Sp(A[\ell], e_\ell^\lambda)$
or $|\Mon_{K'}(A[\ell])|\le 2$. Let $\Lambda$ be the set of prime numbers $\ell\ge \ell_0$
where $|\Mon_{K'}(A[\ell])|\le 2$. We claim that $\Lambda$ is {\em finite}.

For every $\ell\in \Lambda$ we have $[K'(A[\ell]):K']\le 2$. Furthermore
$G(K'(A_{\neq p})/K')$ is profinitely generated, where $p=\chara(K)$. To see this note
that $$G(K'(A_{\neq p})/K')=G(\widetilde{F}K'(A_{\neq p})/\widetilde{F}K')$$ because
$\widetilde{F}/F\sep$ is purely inseparable and use Corollary \ref{fgl}.
Hence $N:=\prod_{\ell\in \Lambda} K'(A[\ell])$ is a {\em finite} extension of $K'$. In
particular $N/F\sep$ is a finitely generated regular extension.
$A/K$ must be geometrically simple by our assumption that $A/K$ has big monodromy
(cf. Proposition \ref{monpropnew1}).
In particular $A_N$ is simple.
Hence assumption that $A$ is not weakly isotrivial with respect to $F$ implies that the
Chow trace ${\mathrm{Tr}}_{N/F\sep}(A_N)$ is zero. It follows by the Mordell-Lang-N\'eron
theorem (cf. \cite[Theorem 2.1]{conrad2006}) that $A(N)$ is a finitely generated
$\Zz$-module.
In particular the torsion group $A(N)\tor$ is finite. On the other hand, $A(N)$ contains a non-trivial
$\ell$-torsion point for every $\ell\in \Lambda$. It follows that $\Lambda$ is in fact finite.

Thus, after replacing $\ell_0$ by a bigger constant, we see that $\Mon_{K'}(A[\ell])=\Sp(A[\ell], e_\ell^\lambda)$
for all primes $\ell\ge \ell_0$. \hfill $\Box$

\begin{coro}\label{ehk}
Let $F$ be a finitely generated field and
$K=F(t)$ the function field of $\Pp^1/F$. Let $A/K$ be a polarized abelian
variety. Let $U\subset \Pp^1$ be an open subscheme such that
$A$ extends to an abelian scheme $\AA/U$. For $u\in U(F)$ denote
by $A_u/F$ the corresponding special fiber of $\AA$.
Assume that $A$ is not weakly isotrivial with respect to $F$ and that
either condition i) or ii) is satisfied.
\begin{enumerate}
\item[i)] $A$ is of Hall type.
\item[ii)] $\chara(K)=0$, $\End(A)=\Zz$ and $\dim(A)=2, 6$ or odd.
\end{enumerate}
Then the sets:
$$X_1:=\{u\in U(F)\,|\quad\End(A_u)\neq \Zz\}$$
and
$$X_2:= \{u\in U(F)\,|\quad\mbox{$A_u/F$ is not geometrically simple}\}$$
are finite.
\end{coro}

{\em Proof.} The abelian variety $A/K$ has big monodromy. In case i) this follows by 
Theorem \ref{main1}. In case ii) this is a well-known theorem of Serre, cf. \cite{serre1984}, 
\cite{serre1985}.) Define $K':=F\sep K$.
As $A/K$ is not weakly isotrivial with respect to $F$ by assumption,
Proposition \ref{geo} implies that there is a constant $\ell_0>\chara(K)$ such
that $\Mon_{K'}(A[\ell])=\Sp(A[\ell], e_\ell^\lambda)$ for all primes $\ell\ge \ell_0$. Hence
$A_{K'}/K'$ has big monodromy. Now Propositions 4 and 7 of \cite{ehk2009}
imply the assertion. Note that the notion of ``big monodromy'' in
the paper \cite{ehk2009} is slightly different from ours.\hfill $\Box$

\section{Appendix. Proof of Theorem \ref{stronghall}}
\label{appendixb}

The aim of this Appendix is to provide a selfcontained proof of Theorem 
\ref{stronghall}, which was first proven in the papers \cite{hall2007} and \cite{hall2008}. We have also taken advantage of the exposition in \cite{Kowalski}.

Let $\ell>2$ be a prime number, let $(V, e)$ be a finite-dimensional symplectic space over $\F_{\ell}$ and $\Gamma=\GSp(V, e)$. In what follows $M$ will be a subgroup of $\Gamma$ which contains a transvection, such that $V$ is a simple $\F_{\ell}[M]$-module. 

\begin{rema}
\begin{itemize}
\item For a set $U\subset V$, we will denote by $\langle U\rangle$ the vector space generated by $U$ in $V$. 
 
\item For a vector $u\in V$ and a scalar $\lambda\in \F_{\ell}$, we denote by $T_u[\lambda]\in \Gamma$ the morphism $v\mapsto v + \lambda e(v, u)u$. For each transvection $\tau\in \Gamma$ there exist $u\not=0$, $\lambda\not=0$ such that $\tau=T_u[\lambda]$, and $u=\ker(\tau-\Id)$. If this is the case we will say that $\langle u\rangle$ is the direction of $\tau$. Each nonzero vector in $\langle u\rangle$ shall be called a \emph{direction vector} of $\tau$.

\item Given a group $G\subset \Gamma$, we will denote by $L(G)$ the set of vectors $u\in V$ such that there exists a transvection in $G$ with direction vector $u$. 

\item We will say that a group $G\subset \Gamma$ fixes a vector space $W$ if $\{g(w): g\in G, w\in W\}\subset W$.
\end{itemize}

\end{rema}

The proof of Part iii) of Theorem \ref{stronghall} is quite simple and is based at the following observation.

\begin{lemm}\label{Partiii}
 Let $G\subseteq \GSp(V)$ be a subgroup and $R$ the subgroup of $G$ generated by the transvections in $G$. Then for all $g\in G$, $r\in R$, $grg^{-1}\in R$.
\end{lemm}

{\em Proof.} Note that if $T=T_v[\lambda]\in G$ is a transvection, then $gT_v[\lambda]g^{-1}=T_{gv}[\lambda]$ is also a transvection, which belongs to $G$, therefore also to $R$. Now if we have an element of $R$, say $T_1\circ \cdots \circ T_k$ for certain transvections $T_1, \dots, T_k$, then $g(T_1\circ \cdots \circ T_k)g^{-1}=(gT_1g^{-1})\circ \cdots \circ  (gT_kg^{-1})$ is the composition of transvections of $G$, therefore an element of $R$.
\hfill $\Box$

Part i) of Theorem \ref{stronghall} is essentially Lemma 3.2 of \cite{hall2007}. Before proceeding to prove it, note the following elementary facts.

\begin{lemm}\label{LemaSuma}
 Let $G$ be a group that acts irreducibly on $V$, and let $W\subset V$ a nonzero vector space. Then $V=\sum_{g\in G} gW$.
\end{lemm}

{\em Proof.}  Let $S$ be the set $S=\{g(w): g\in G, w\in W\}$. Consider the vector space $\langle S\rangle$. This vector space is fixed by $G$, hence since $G$ acts irreducibly on $V$ it must coincide with $V$. 
\hfill $\Box$

\begin{lemm}\label{DivisionOfTransvections}
 Let $W$ be a vector subspace of $V$, and assume that it is fixed by a transvection $T=T_u[\lambda]$. Then either $u\in W$ or $u\in W^{\perp}$.
\end{lemm}

{\em Proof.}
Recall that, for all $v\in V$, $T(v)=v +\lambda e(v, u)u$. If $u\not\in W$, the only way for $T$ to fix $W$ is that $e(w, u)=0$ for all $w\in W$. 
\hfill $\Box$

{\em Proof of Theorem \ref{stronghall}, i)}

Consider the action of $R$ on $V$. The first step is to fix one simple nonzero $R$-submodule $W$ 
contained in $V$ (This always exists because $V$ is finite-dimensional as an $\F_{\ell}$-vector space).

By Lemma \ref{LemaSuma}, we know that $V=\sum_{g\in M} gW$. Moreover, for $g_1, g_2\in M$ it holds  that $g_1W=g_2W$ if and only if  $g_1H= g_2H$. Therefore we can write $V=\sum_{g\in M/H} gW$, where $H$ is the stabilizer of $W$ in $M$. The proof of i) boils down to prove that the sum is direct and orthogonal, that is, if $g_1H\not=g_2H$, then $g_1W\cap g_2W=0$ and $g_1W\subset (g_2W)^{\perp}$. Equivalently, we will prove that for any $g\in M$, if $gW\not=W$, then $gW\cap W=0$ and $gW\perp W$.

The first claim, namely $gW\not=W$ implies $gW\cap W=0$ is easy. The key point is to note that for each $g\in M$, $gW$ is also fixed by $R$. Take $r\in R$, $gw\in gW$. Then $rgw=g(g^{-1}rg)w\in gW$ since $g^{-1}rg\in R$ by Lemma \ref{Partiii} and hence fixes $W$. Now it follows that $W\cap gW$ is fixed by $R$, and thus is an $R$-subrepresentation of $W$. But $W$ is an simple $R$-module, hence since $W\cap gW\not=W$, it must follow that $gW\cap W=0$. 

To prove that $gW\not=W$ implies $gW\perp W$, we need to make first the following very important observation.

\begin{claim}\label{Claim1}
 The set $L(M)\cap W$ generates $W$.
\end{claim}

{\em Proof of Claim \ref{Claim1}.}
First let us see that $L(M)\cap W$ is nontrivial. Since any  transvection in $M$ fixes $W$ by definition of $W$, it follows by Lemma \ref{DivisionOfTransvections} that either its direction vector belongs to $W$, or else it is orthogonal to $W$, in which case the transvection acts trivially on $W$. But it cannot happen that all transvections in $M$ act trivially on $W$. For, if a transvection $T$ acts trivially on $W$, then for all $g\in M$, $gTg^{-1}$ acts trivially on $gW$. But since $R=gRg^{-1}$  (because of Lemma \ref{Partiii}), then if all $R$ acts trivially on $W$, it also acts trivially on $gW$. Now recall that $V=\sum_{g\in M}gW$. Then $R$ would act trivially on $V$. But $R$ contains at least a  transvection, and this does not act trivially on $V$. We have a contradiction.

Hence $L(M)\cap W$ is non zero. But now observe that this set is fixed by the action of $R$, since the elements of $M$ bring direction vectors into direction vectors. Therefore the vector space $\langle L(M)\cap W\rangle\subset W$ is fixed by the action of $R$. Since we are assuming $W$ is an simple $R$-module, it follows that $\langle L(M)\cap W\rangle= W$
\hfill $\Box$

Now we are able to prove that if $gW\not=W$, then $gW\subset W^{\perp}$.  Because of the previous claim, it suffices to show that, for any nonzero vector $w\in W$ which is the direction vector of a transvection in $M$, say $T$, $w\in (gW)^{\perp}$. Now recall that, since $T$ fixes $gW$, by Lemma \ref{DivisionOfTransvections} either $w\in gW$ or $w\in (gW)^{\perp}$. But $gW\cap W=0$, so $w\in (gW)^{\perp}$.
\hfill $\Box$

Before proving Part ii) of Theorem \ref{stronghall}, we will introduce some notation.

\begin{defi}
 Let $g\in M$. We will denote by $R_g$ the subgroup of $R$ generated by the transvections that act non-trivially on $gW$.
\end{defi}

The following lemma is Lemma 7 of \cite{hall2008}.

\begin{lemm}\label{Commutativity} Let $g_1, g_2\in M$ with $g_1H\not=g_2H$. Then the commutator $[R_{g_1}, R_{g_2}]$ is trivial.
\end{lemm}

{\em Proof.}
For $i=1, 2$, let $T_i\in R_{g_i}$ be a transvection. We will see that they commute. By Lemma \ref{DivisionOfTransvections} applied to $g_iW$, either $T_i$ acts trivially on $g_iW$ or its direction vector, say $u_i$, belongs to $g_iW$. By definition of $R_{g_i}$ we have the second possibility. But because of Part i) of Theorem \ref{stronghall}, for each $g\in M$ such that $g_iW\not=gW$, $g_iW\cap gW=0$, hence $u_i\not\in gW$. Therefore again by Lemma \ref{DivisionOfTransvections} applied now to $gW$, it follows that $T_i$ acts trivially on $gW$. Therefore $T_1$ and $T_2$ commute on each $gW$, since at least one of them acts trivially on it. Since $V=\bigoplus_{g\in M/H} gW$, it follows that they commute on all $V$.    
\hfill $\Box$

{\em Proof of Theorem \ref{stronghall}, ii).}
 
Let $M/H=\{g_1H, \dots, g_sH\}$, with $g_1=\Id$. Define the map
\begin{equation*}\begin{aligned} P:\prod_{i=1}^s R_{g_i} &\rightarrow R\\
                                   (r_1, r_2, \dots, r_s)&\mapsto r_1\cdot r_2\cdot\cdots\cdot r_s.\end{aligned}\end{equation*}

Since by Lemma \ref{Commutativity} elements from the different $R_{g_i}$ commute, this map is a group homomorphism. Let us see that it is also an isomorphism.

Assume that $r_1\cdot r_2\cdot \cdots r_s=\Id$, and that there is a certain $r_j$ which is not the identity matrix. Then $r_j$ must act nontrivially on a certain vector $v\in V$. Since the elements of $R_{g_j}$ act trivially on the elements of $g_iW$ for $i\not=j$ and $V=\bigoplus_{i=1}^s g_iW$, we can assume that $v\in g_jW$. But then the remaining $r_i$ with $i\not=j$ act trivially on $v$ and on $r_j(v)$. Therefore $\Id(v)=r_1\cdot \cdots \cdot r_s(v)=r_j(v)\not=v$, which is a contradiction. To prove surjectivity, it suffices to note that each transvection $T$ of $M$ belongs to one of the $R_{g_i}$, (hence each element of $R$ can be generated by elements of $\cup_i R_{g_i}$). And this holds because, since $T$ fixes all the $g_iW$, the direction vector of $T$ must either belong to $g_iW$ or be orthogonal to it because of Lemma \ref{DivisionOfTransvections}, and since $V=\oplus_{i=1}^s g_iW$ it cannot be orthogonal to all the $g_iW$. Therefore we get that $R\simeq \prod_{i=1}^s R_{g_i}$.

Now we are going to apply the following result \cite[Main Theorem]{ZS1}:

\begin{thm}
Suppose $G \subset \GL(n, k)$ is an irreducible group generated by
transvections. Suppose also that $k$ is a finite field of
characteristic $\ell > 2$, and that $n>2$. Then $G$ is conjugate in
$\GL(n, k)$ to one of the groups $\SL(n,k_0)$, $\Sp(n,k_0)$ or
$\mathrm{SU}(n,k_0)$, where $k_0$ is a subfield of $k$.
\end{thm}

Note that, if $n=2$, the result is also true and well known (cf. \cite[Section 252]{Dickson}). 

Now $R_{g_1}$ is generated by transvections, and acts irreducibly on $W$ (because $R$ acts irreducibly on $W$, and $R_{g_1}$ is the group generated by all those transvections in $M$ that act nontrivially on $W$). Therefore $R_{g_1}$ is conjugated to $\Sp(W)$. Since all $R_{g_i}$ are conjugated to $R_{g_1}$, the same holds for them. Therefore we have the isomorphism
$R\simeq \prod_{i=1}^s \Sp(W)$.

Finally, we can view $H_1=\prod_{i=1}^s \GSp(W)\simeq\prod_{i=1}^s \GSp(g_iW)$ as the subgroup of $\Gamma$ fixing each $g_iW$ and, fixing a symplectic basis on each $g_iW$, we can view $H_2=\Sym(M/H)$ as the subgroup of $\Gamma$ that permutes the $g_iW$ by bringing the fixed symplectic basis of each $g_iW$ into the fixed symplectic basis of another $g_jW$. The group generated by $H_1$ and $H_2$ inside $\Gamma$, which is the group of elements of $\Gamma$ that permute the $g_iW$, is the semidirect product $H_1\rtimes H_2$.

Recall that $N_{\Gamma}(R)=\{g\in \Gamma: g R g^{-1}=R\}$. Note that $g\in N_{\Gamma}(R)$ if and only if for all transvections $T\in M$, $gTg^{-1}\in R$. Now, if $T=T_v[\lambda]$, it holds that $gTg^{-1}=T_{g(v)}[\lambda]$, and this transvection belongs to $R$ if and only if it is a transvection of $M$, that is to say, if and only if $g(v)\in L(M)$. Therefore $g\in N_{\Gamma}(R)$ if and only if $g(L(M))=L(M)$. Now since $R$ is isomorphic to $\prod_{i=1}^s\Sp(g_iW)$, $L(M)$ is the disjoint union of the $g_iW$. And moreover, if $W$ is an  $R$-module and $g\in N_{\Gamma}(R)$, then $R$ fixes $gW$. Therefore, if $W$ is an simple $R$-module, then 
$gW\not=W$ implies that $gW\cap W=0$. Thus if $g\in N_{\Gamma}(R)$, then $g$ permutes the $g_iW$. In other words, $N_{\Gamma}(R)\subset \prod_{i=1}^s \GSp(W)\rtimes \Sym(M/H)$. Reciprocally, each element of $\prod_{i=1}^s \GSp(W)\rtimes \Sym(G/H)$ carries elements of $\bigcup_i g_iW$ in elements of $\bigcup_i g_iW$, that is to say, carries $L(M)$ into $L(M)$, and therefore belongs to $N_{\Gamma}(R)$.
\hfill $\Box$

This completes the proof of Part a) of Theorem \ref{stronghall}.

{\em Proof of Part b) of Theorem \ref{stronghall}.} Recall that $(V, e)$ is a symplectic space over
$\Ff_\ell$ and $M$ a subgroup of $\Gamma:=\GSp(V, e)$. $M$ contains
a transvection and $V$ is a simple $\Ff_\ell[M]$-module by assumption. Furthermore 
$R$ is the subgroup of $M$ generated by the transvections in $M$, 
$0\neq W\subset V$ is
a simple $\Ff_\ell[R]$-module and $H=\Stab_M(W)$. We already proved that
there is a orthogonal direct
sum decomposition $V=\bigoplus_{g\in M/H} gW$. Furthermore 
$R\cong \prod_{g\in M/H} \Sp(W)$, $N_\Gamma(R)\cong \prod_{g\in M/H} \GSp(W)\rtimes
\Sym(M/H)$ and $R\subset M\subset N_\Gamma(R)$.
Denote by $\varphi: N_\Gamma(R)\to \Sym(M/H)$ the projection.

Let $E/\Ff_\ell$ be a finite extension and $\rho: E^\times\to M\subset \GL(V)$ a 
representation of amplitude $\amp(\rho)\le e$. Assume that $\ell> e\dim(V)+1$.
We have to prove that $\varphi(\rho(E^\times))=\{1\}$. 

Define $S:=\ke(\varphi\circ \rho)\subset E^\times$. Then $[E^\times:S]\le |M/H|\le \dim(V)$, and
this implies $e[E^\times:S]<\ell-1$.
Furthermore 
$$\rho(S)\subset \ke(\varphi)\cong\prod_{g\in M/H} \GSp(gW).$$
Obviously $\rho(S)$ commutes with the center
$$Z(\ke(\rho))\cong\prod_{g\in M/H} \Ff_\ell^\times Id_{gW}$$
of $\ke(\rho)$. Now by \cite[Lemma 3]{hall2008} $\rho(E^\times)$ commutes with
$Z(\ke(\rho))$, because $e[E^\times:S]<\ell-1$. It can easily be seen that the 
centralizer of $Z(\ke(\rho))$ in $N_\Gamma(R)$ is equal to 
$\ke(\varphi)\cong \prod_{g\in M/H} \GSp(gW)$. Hence $\rho(E^\times)\subset \ke(\varphi)$ and
this implies $\varphi\circ \rho(E^\times)=\{1\}$.\hfill $\Box$
\bigskip\bigskip

\noindent
{\bf Acknowledgements.}
S. A. is a research fellow of the Alexander von Humboldt Foundation.
S. A. was partially supported by the Ministerio de Educaci\'on y Ciencia grant MTM2009-07024. S. A. wants to thank the Hausdorff Research Institute for Mathematics in Bonn, the Centre de Recerca Matem\`atica in Bellaterra and the Mathematics Department of Adam Mickiewicz University in Pozna\'n for their support and hospitality while she worked on this project.
W.G. and S.P. were partially supported by the Deutsche Forschungsgemeinschaft
research grant GR 998/5-1.
W.G. was partially supported by the Alexander von Humboldt Research Fellowship
and an MNiSzW grant. W.G. thanks Centre Recerca Matem\`atica in Bellaterra and the Max Planck
Institut f{\" u}r Mathematik in Bonn for support and hospitality during visits in 2010, when he worked on this
project. S.P. gratefully acknowledges the hospitality of Mathematics Department of Adam Mickiewicz  University in Pozna{\' n} and of the Minkowski center at Tel Aviv University during several research visits.

\bibliographystyle{plain}
\bibliography{lit}

\end{document}